\documentclass[11pt,leqno]{article}

\usepackage{graphicx}

\newcommand{\maths}[1]{{\bf #1}}

\newtheorem{defi}{Definition}[section]
\newtheorem{lemm}[defi]{Lemma}
\newtheorem{prop}[defi]{Proposition}
\newtheorem{theo}[defi]{Theorem}
\newtheorem{coro}[defi]{Corollary}
\newtheorem{conj}[defi]{Conjecture}

\newcommand{\bdefi}{\begin{defi}}
\newcommand{\edefi}{\end{defi}}
\newcommand{\bprop}{\begin{prop}}
\newcommand{\eprop}{\end{prop}}
\newcommand{\btheo}{\begin{theo}}
\newcommand{\etheo}{\end{theo}}
\newcommand{\blemm}{\begin{lemm}}
\newcommand{\elemm}{\end{lemm}}
\newcommand{\brema}{\begin{rema}}
\newcommand{\erema}{\end{rema}}
\newcommand{\bexer}{\begin{exam}}
\newcommand{\eexer}{\end{exam}}
\newcommand{\bconj}{\begin{conj}}
\newcommand{\econj}{\end{conj}}
\newcommand{\bcoro}{\begin{coro}}
\newcommand{\ecoro}{\end{coro}}

\newcommand{\dem}{\noindent{\bf Proof. }}
\newcommand{\rem}{\noindent{\bf Remark. }}

\newcommand{\rems}{\noindent{\bf Remarks. }}
\newcommand{\que}{\noindent{\bf Question. }}

\newcommand{\T}{{\cal T}}

\newcommand{\ra}{\rightarrow}

\newcommand{\RR}{\maths{R}}

\newcommand{\bord}{\partial}
\newcommand{\df}{d\!f}
\newcommand{\du}{d\!u}

\newcounter{fig}

\def\myfigurepdf#1#2#3{
\addtocounter{fig}{1}
\[
\begin{array}{c}
\mbox{\includegraphics[width=#2]{#1}}\\
\\
\hbox{\rm Figure \arabic{fig} ~: #3.}
\end{array}
\]
}

\def\myfigurepdft#1#2#3{
\addtocounter{fig}{1}
\[
\begin{array}{c}
\mbox{\includegraphics[width=#2]{#1}}\\
\\
\hbox{\rm Figure \arabic{fig} ~ #3}
\end{array}
\]
}

\newcommand{\eop}[1]{{\flushright\hfill\fbox{\bf #1}}}
\newcommand{\ack}{\noindent{\bf Acknowledgement.}}

\title{Energy and length in a topological planar quadrilateral}

\author{Sa'ar Hersonsky
\thanks{The research is partially supported by
The David \& Elaine Potter Charitable Foundation}}

\date{}

\begin{document}

\maketitle

\begin{abstract}
\noindent  We provide bounds for the product
of the lengths of 
distinguished shortest paths 
 in a finite network induced by a triangulation of a topological planar quadrilateral.

\end{abstract}

\section{Introduction}
\label{sec:Intro}

A topological planar closed disk with four
distinguished points on its boundary, its corners, will be called a quadrilateral. The following definition is due to Schramm (\cite{Sch}).
\bdefi 
\label{de:triangulated quad}{\rm Let $Q$ be a 
quadrilateral endowed with a triangulation. Let $V,E,T$ denote the set of vertices, edges and
triangles of $Q$, respectively. Let $\bord Q=P_1\bigcup P_2\bigcup
P_3\bigcup P_4$ be a decomposition of $\bord Q$ into four 
non-trivial arcs of the triangulation with disjoint interiors, in cyclic order. If the intersection of two any of these arcs is not empty, then it consists of a corner (all of which are vertices).
 A corner must belong to one and
only one of the $P_i$'s.  
The collection
$\T=(V,E,T,P_1,P_2,P_3,P_4)$ will be called a triangulation of
$Q$. }
\edefi 
By invoking a {\it conductance function}, $\T$ becomes a {\it finite network}. One can define
a {\it boundary value problem} (BVP) on the network. 
  Let $f$ be the 
solution of the BVP and let $I(f)$ be its {\it Dirichlet energy}. Corollary~\ref{co:product of lengths} provides inequalities relating the product of the lengths of a shortest {\it thick} vertical path (a particular path which connects $P_3$ and $P_1$) and a shortest thick horizontal path (a particular path which connects $P_2$ and $P_4$)  in terms of  $I(f)$ and some constants arising from the combinatorics and the conductance function. 
 Corollary~\ref{co:product of lengths} follows from Theorem~\ref{le:length estimate} (our main theorem) and Lemma~\ref{le:upper bound}. 
 The length is measured with respect to
 $\rho$,  the {\it gradient metric}  (Defintion~\ref{de:gradient metric}) which is induced by the solution of the BVP (see Section~\ref{se:Preliminaries} and Section~\ref{se:length-energy estimates} for the precise definitions of the notions above).
In the special case where $c(x,y)\equiv 1$ and $k$ is the maximal degree of $V$, it follows from Corollary~\ref{co:product of lengths} that 
$$ l(|V|) I(f)\geq {\mbox{\rm length}}_{\rho}(\gamma^*){\mbox{\rm length}}_{\rho}(\gamma)\geq \frac{1}{\sqrt k}I(f),\   \mbox{\rm  where}\  l(|V|)\  \mbox{\rm is some constant}.  $$

\medskip
C.~Loewner (see \cite{Bes}) studied  Differential-geometric inequalities relating area and the product of shortest (vertical and horizontal) curves in a quadrilateral. His inequalities are derived with respect to the Euclidean metric.
His work was generalized and forms a rich theory. The well known reciprocal property of the {\it extremal lengths} of conjugate families of curves in a quadrilateral is one useful example of this theory (see \cite{AhSa},  \cite{MaRo} and \cite{Oh} for a few examples and generalizations for other Riemann surfaces). 
 
 The original notion of extremal length in a discrete setting was introduced by Duffin (\cite{Du}). More recently, in this setting, Cannon  (\cite{Ca}) introduced a different notion of extremal length.
In the work of Cannon, Floyd and Parry (see for instance \cite{CaFlPa}), as well as in the work of Schramm (\cite{Sch}), inequalities generalizing the reciprocal property of the extremal length (with respect to some extremal metric) of conjugate families of curves in a quadrilateral are very useful.

In the setting of finite and infinite networks reciprocal property of extremal length and {\it capacity} were studied extensively with respect to an extremal metric (see \cite{So} for a detailed account).

One motivation for using the gradient metric in this paper arises from  extremal length arguments in the complex plane. It is well known (see \cite{AhSa}) that for every $z\in {\bf C}$   the extremal metric in a topological quadrilateral in the complex plane satisfies $$m(z)=| \nabla(f)(z)|,$$ where $f$ is the solution of the classical Dirichlet-Neumann boundary value problem. In the complex plane it is also known that equality holds in Corollary~\ref{co:product of lengths}, where both sides equal $I(f)$.

Consider the broader class of  BVP problems that are studied in \cite{BeCaEn}.  Our work is also motivated by the following. \newline
\que
Is there a BVP problem and a metric $\rho_0$ (which is perhaps different than $\rho$) derived from the solution such that 
$$ {\mbox{\rm length}}_{\rho_0}(\gamma^*){\mbox{\rm length}}_{\rho_0}(\gamma)\geq I(f)\ \mbox{\rm and} \sum_{x\in\bar F} \rho_{0}^2(x)=I(f)     \mbox{\rm (see Definition~\ref{de:discrete metric}}) ?$$ 

\rem
In paper \cite{Her} we will use some of the ideas of this paper to prove a
finite Riemann mapping theorem (\cite{CaFlPa}, \cite{Sch}). A more direct proof would follow 
from a positive answer to the question above.
A finite Riemann mapping theorem can be viewed as the first step in solving the Cannon conjecture: A negatively curved group $G$ with $\partial G=S^2$ is Kleinian. We hope that our ideas will be useful towards the resolution of this conjecture.

\medskip

\noindent{\small \ack  \newline Part of this research was done
while the author visited the IAS, the  Mathematics departments of Princeton
and Caltech during the summer of 2003. We express our deepest
gratitude to both for their generous hospitality. We
are grateful to Francis Bonahon,
 Nikolai Makarov and Steven Kerckhoff for enjoyable and helpful discussions. We are indebted to David Gabai and Fr\'ed\'eric Paulin for their careful reading of this work and many useful remarks. We are indebted to anonymous referees of this paper. Their remarks and corrections were essential. 
}

\section{Preliminaries}
\label{se:Preliminaries}

We recall some known facts regarding harmonic functions
 and boundary value problems on networks.
We use the notations of Section 2 in \cite{BeCaEn}. Let
$\Gamma=(V,E,c)$ be a {\it finite network}, that is a simple and
finite connected graph with a vertex set $V$ and edge set $E$.
We shall also assume that the graph is planar.
 Each edge $(x,y)\in E$ is
assigned a {\it conductance} $c(x,y)=c(y,x)>0$.
 Let ${\cal P}({ V})$
denote the set of non-negative functions on $V$. If $u\in {\cal
P}( V)$, its support is given by $S(u)=\{ x \in V: u(x)\neq 0 \}$.
Given $F\subset V$ we denote by $F^{c}$ its complement in
$V$.  Set
${\cal P}(F)=\{u\in {\cal
P}(V):S(u)\subset F\}$.  The set   $\bord F=\{ (x,y)\in E: x\in F,  y\in F^{c} \}$  is called 
the {\it
edge boundary} of $F$ and the set  $\delta F=\{x\in F^{c}: (x,y)\in E\ {\mbox
{\rm for some}}\ y\in F \}$ is called the {\it vertex boundary} of
$F$. Let ${\bar F}=F\bigcup \delta F$ and let $\bar E=\{(x,y)\in
E :x\in F\}$.
Given $F \subset V$, let
${\bar \Gamma}(F)=({\bar F},{\bar E},{\bar c})$ be the network
such that ${\bar c}$ is the restriction of $c$ to ${\bar E}$. 
We say that $x\sim y$ if $(x,y)\in \bar E$. For $x\in \bar F$ let $k(x)$ denote the degree of $x$ (if $x\in\delta(F)$ the neighbors of $x$ are taken
  only  from $F$).

For $f,h:\bar E \rightarrow R$ we let 
 $(f,h)=\sum_{(x,y)\in \bar E}\frac{f(x,y)h(x,y)}{c(x,y)}$ be an
inner product on $l^2(\bar E,1/c)$ (see \cite[1.2.A]{Woe}).
The following definitions are discrete analogues of
classical notions in continuous Potential Theory \cite{Fu}.


\bdefi \mbox{ \rm (\bf\cite[Section 3]{BeCaEn1})}
{\rm \label{def:energy}  Let $u\in {\cal P}({\bar  F})$, 
\begin{enumerate}
\item then for $x\in \bar F$, the function $\Delta u(x)=\sum_{y\sim x}c(x,y)\left( u(x)-u(y) \right )$ is called
  the potential of $u$ at $x$, \mbox{\rm(}if $x\in\delta(F)$ the neighbors of $x$ are taken
  only  from $F$\mbox{\rm)} and
\item the number
$I(u)= \sum_{x\in\bar F}\Delta u(x)u(x)=\sum_{(x,y)\in \bar E} c(x,y)( u(x)-u(y)
)^2,$ is called the {\it Dirichlet energy} of $u$.

\item A function $u\in {\cal P}({\bar F})$ is called harmonic in $F\subset V$ if
$\Delta u(x)=0,$ for all $x\in F$.
\end{enumerate}
}
\edefi

When $c(x,y)\equiv 1$, an easy computation
shows that $u$ is harmonic at a vertex $x$ if and only if the
value of $u$ at $x$ is the arithmetic average of the value of $u$
on the neighbors of $x$. 

When $(x,y)\in \bar E$ let us denote by $[x,y]$ the directed
edge from $x$ to $y$ and let $\overrightarrow E=\{[x,y] : (x,y)\in \bar E \}$ denote the  set of all directed
edges. Given $u : V\ra \RR$ we define the {\it differential } or
the {\it gradient} of $u$ as $\du: \overrightarrow E\ra \RR$ by
$\du[x,y]=c(x,y)(u(y)-u(x))$ for all $[x,y]\in \overrightarrow E$ (see for instance the notations of Section 2 in \cite{BeSch2}).
Note that if
$|\overrightarrow E|=m$, then $\du$ can be identified with a vector in $\RR^m$.
 It now follows by
Definition~\ref{def:energy} that for every function
$u:V\ra \RR$ we have that 
$I(u)= \frac{1}{2}\sum_{e \in \overrightarrow E}\|\du( e)\|^2.$

Let $\overrightarrow E(x)$ denote the set of all
edges of the form $[x,y]$ which are in $\overrightarrow E$. Any
$g:\overrightarrow E(x)\ra \RR$ can be naturally viewed as an element in
$\RR^{k(x)}$. We will denote this vector space, with the restriction of the inner
product on $\bar E$, by $T_x$. In particular we have

\bdefi
{\rm 
\label{de:restricted gradient} Let $u:\bar F\ra \RR$.  Let
$ \overrightarrow \du(x)=\du|_{\overrightarrow E(x)}$, denote the restriction of $\du$ to
$\overrightarrow E(x)$ (In particular, $\overrightarrow\du(x)$ can be viewed as a vector in
$T_x $).}
\edefi

For $x\in \delta(F)$ let $\{y_1,y_2,\ldots,y_m\}\in F$ be its neighbors, enumerated in a cyclic order.
\bdefi
{\rm 
\label{de:vector normal derivative} 
The normal vector derivative at $x\in \delta(F)$ is defined by \newline
$\overrightarrow{\frac{\bord u}{\bord n_{F}}}(x)=
\left( c(x,y_1)(u(x)-u(y_1)),\ldots ,c(x,y_m)(u(x)-u(y_m) ) \right )$
  and the conductance vector at $x$ is defined by
$\overrightarrow{c}_{\delta(F)}(x)=\left(c(x,y_1),\ldots,c(x,y_m)\right)$. 
}
\edefi
 If $x\in F$, $\overrightarrow {c}_{F}(x)$ is defined similarly and the neighbors of $x$ are taken in $F \bigcup \delta(F)$.

The following definition provides the discrete
analogue of the continuous notion of {\it normal derivative}.

\bdefi \mbox{\rm \bf(\cite{ChGrYa})}
{\rm 
\label{def:normal}
The normal derivative of $u$ at a point
$x\in \delta F$ with respect to the set $F$ is $$\frac{\bord u}{\bord n_{F}}(x)= \sum_{y\sim x,\
y\in F}c(x,y)  (u(x)-u(y) ) .$$ 
}
\edefi

\smallskip 

The following
 proposition establishes a discrete version of the first classical {\it Green
identity}. It will be crucial in the proof of  Theorem~\ref{le:length estimate}.
\bprop \mbox{\rm \bf(\cite[Prop. 3.1]{BeCaEn}) (The first Green
identity)}
\label{pr:Green id} Let $F \subset V$ and $u,v\in {\cal P}({\bar
F})$. Then we have that
 $$ \sum_{(x,y)\in {\bar
E}}c(x,y)(u(x)-u(y))(v(x)-v(y))=\sum_{x\in F}\Delta
u(x)v(x)+\sum_{x\in\delta(F)}\frac{\bord u}{\bord n_{F}}(x)v(x).$$
\eprop

\rems

\noindent  1. In \cite{BeCaEn} a second Green identity is obtained. In this paper we will
   use only the one above.

\noindent 2. In \cite{BeCaEn3} (see in particular Section 2 and Section 3) a systematic study of 
discrete calculus on $n$-dimensional (uniform) grids of Euclidean $n$-space is provided. Their definition of a tangent space may be adopted to our setting and does not require the notion of directed edges. However, in \cite{Her} directed edges will play an important role.

\section{Length estimates of shortest paths}
\label{se:length-energy estimates}
Throughout this section  $\T$ 
 will denote a fixed triangulation of a quadrilateral
(see Definition~\ref{de:triangulated quad}). We will denote by $F$  the set of
vertices which do not belong to $\bord Q$. Hence, $\delta(F)$ is
the set of vertices that belong to $P_1\bigcup P_2\bigcup P_3
\bigcup P_4$. Let  $\{c(x,y)\}_{(x,y)\in \bar{E}}$ be a fixed conductance function and let $\bar{\Gamma}( F)$ be the associated network.
We are interested in functions that solve a boundary value problem (BVP) on $\bar{\Gamma}( F)$. The following definition is based on \cite[Section 3]{BeCaEn} and \cite[Section 4]{BeCaEn2}.
\bdefi
{\rm 
\label{de:boundary function}
 Let $g>0$ be a constant.
 A Dirichlet-Neumann boundary value function
 is a  function $f \in {\cal P}({\bar F})$ which satisfies the following:
\begin{enumerate}
  \item $f$ is harmonic in $F$,
  \item $f|_{P_{2}}=0$,
  \item $f|_{P_{4}}=g$, for some constant $g$,  and
  \item $\frac{\bord f}{\bord n_{F}}|_{P_{1}}=\frac{\bord f}{\bord n_{F}}|_{P_{3}}=0$.
\end{enumerate}
}
\edefi
\rem
The uniqueness and existence of a Dirichlet-Neumann boundary value
function is provided by the nice and foundational work in \cite[Section 3]{BeCaEn} and  \cite[Section 4]{BeCaEn2}. In fact, their work provides a detailed framework for a broader class of boundary value problems on finite networks.

\bdefi  \mbox{ \rm (\bf\cite{Ca})}
{\rm 
\label{de:discrete metric}
A metric on a finite network is a function $\rho : V\ra [0,\infty)$.
}
\edefi
In particular, the length of a path is given by integrating $\rho$ along the path. 
When $\rho\equiv 1$,  the familiar distance function on $V\times V$ is obtained by setting 
$\mbox{\rm dist}( A, B)= \sum_{x\in \alpha} 1 -1=  k$, where $\alpha=(x,x_1,\ldots x_k)$ is a path with the smallest possible number 
of vertices among all the paths connecting a vertex in $ A$ and a vertex in $B$. 
We now define the gradient metric which will be used in our estimates. 
\bdefi 
{\rm 
\label{de:gradient metric}
Given $f\in {\cal P}({\bar F})$ the gradient metric induced by $f\in {\cal P}({\bar F})$  is defined by
\label{def:metric}
 $$\rho(x)=\left \{ \begin{array}{ll}
    \|\overrightarrow\df(x)\|&\mbox{\rm if}\  x\in F \\
    \\
    \|\overrightarrow{ \frac {\bord f}{\bord n_{F}}}(x)\|&\mbox{\rm if}\  x\in
    \delta(F).
  \end{array} \right.$$ 
  
  }
\edefi

Before turning to our main theorem, we will define the paths which are going to be considered.  The definition below describes  two classes of paths. These classes are sufficiently separated from $\delta(F)$.   

\bdefi 
{\rm
\label{de:thick path}
A path $\beta=(x_0,x_1,\ldots,x_n)$ will be called vertically  (horizontally) thick if it satisfies the following:
\begin{enumerate}
\item $x_0= \beta \cap P_3$ ($x_0= \beta \cap P_2$) and $x_n= \beta \cap P_1$
 ($x_n= \beta \cap P_4$) respectively,
 
\item neither $x_0$ or $x_n$ is a corner, 

\item for all $i=1,\ldots n-1$, $x_i\in F$, 

\item for all $i=2,\ldots n-2$, $\mbox{\rm dist}( x_i, P_3\cup P_4\cup P_1)> 1$,

\item $\mbox{\rm dist}(x_1, x), \mbox{\rm dist}(x_{n-1},x)\geq 1$ when $x\in P_1\cup P_3\cup P_4$ 
 ($\mbox{\rm dist}(x_1, x),\mbox{\rm dist}(x_{n-1},x) \geq 1$ when $x\in\delta(F)$), and equality is attained uniquely for $x_0$ and $x_n$ respectively. 

\end{enumerate}
}
\edefi

We now turn to our main theorem.

\btheo 
\label{le:length estimate}
Let $\T$ be a triangulation of a topological quadrilateral. Let  $\bar{\Gamma}( F)$ be the associated network.
Let $f$ be the Dirichlet-Neumann boundary value function with some constant $g$.  Let $\rho$ be the gradient metric induced by $f$.
Let 

 $M=\max_{x\in F\cup\delta(F)}\{  \|\overrightarrow{c}_{F}(x)\|, \|\overrightarrow{c}_{\delta(F)}(x)\}\|\}$  and let $m=\min_{(x,y)\in
\bar E}\sqrt {c(x,y)}$.
\begin{enumerate}
\item  If $\gamma^*$ is a shortest vertical thick path which connects $P_1$ to $P_3$, then
$${\mbox{\rm length}}_{\rho}(\gamma^*)\geq \frac{
I(f)}{g M},\  {\mbox{\rm and}}$$
\item if $\gamma$ is a shortest horizontal thick path which connects $P_2$ to $P_4$, then
   $${\mbox{\rm length}}_{\rho}(\gamma)\geq g m .$$
\end{enumerate}
\etheo

\dem Using properties (1)-(4) of $f$ and the
first Green identity (Proposition~\ref{pr:Green id}) with
$u=v=f$ we obtain that
\begin{equation}
\label{eq:properties}
  I(f)=\sum_{x\in P_4}\frac{\bord f}{\bord n_{F}}(x)f(x).
\end{equation}

Hence, by the definition of
$g$,  we have that
\begin{equation}
\label{eq:properties1}
  I(f) =  g\  |\sum_{x\in P_4}\frac{\bord f}{\bord n_{F}}(x)|.
\end{equation}

Let $\gamma^{*}$ be a shortest thick path connecting $P_3$ to $P_1$.
It is clear that we may assume that $\gamma^{*}$ is simple. Let $ x_0= \gamma^{*}\cap P_3$ and  let
$x_n=\gamma^{*}\cap P_1$.

\myfigurepdf{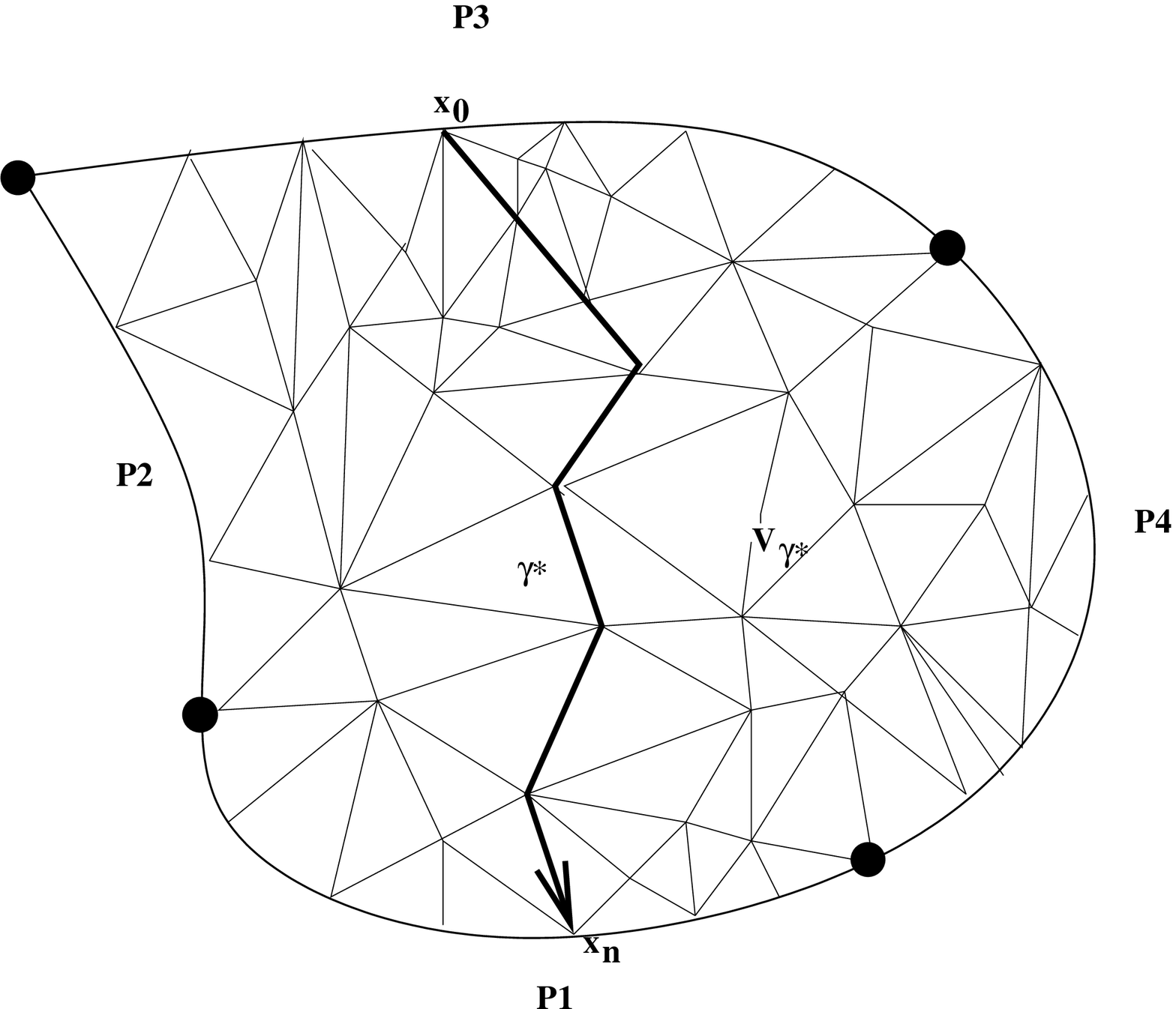}{7.5cm}{A shortest thick vertical path between $P_3$ and $P_1$}

Let $V_{\gamma^{*}}$ denote the subset of $F$ which is
   enclosed, in a cyclic order,  by
$\gamma^{*}$, a part of $P_1$ (which we will denote by $P_{1}(\gamma^{*})$),  $P_4$ and  part of $P_3$ (which we will denote by $P_{3}(\gamma^{*})$).
It follows from Definition~\ref{de:thick path}   that $V_{\gamma^{*}}\neq \emptyset$ and that 
$\delta(V_{\gamma^{*}})=\gamma^{*}\cup  P_{1}(\gamma^{*}) \cup P_4 \cup P_{3}(\gamma^{*})$.

We now apply the first Green identity
 with $u=f$ and
  the constant function $v \equiv 1$ in 
$V_{\gamma^{*}}\cup\delta (V_{\gamma^{*}})$. We obtain that
\begin{equation}\sum_{x\in \delta (V_{\gamma^{*}})} \frac{\bord f}{\bord n_{V_{\gamma^{*}}}}(x)
=0.
\end{equation}

It follows by Definition~\ref{de:thick path}  that for every $x\in P_{1}(\gamma^{*})\cup P_4 \cup P_{3}(\gamma^{*})$ which is not in $\gamma^{*}$ we have that $$\frac{\bord f}{\bord n_{V_{\gamma^{*} } }  }  (x)=\frac{\bord f}{\bord n_{F } }(x)   .$$

By using the fourth property of Definition~\ref{de:boundary function}  and the triangle inequality we have that
\begin{equation}
\label{eq:stokes}
   | \sum_{x\in P_4}\frac{\bord f}{\bord n_{V_{\gamma^{*} } }  }  (x) |=
  |\sum_{x\in\gamma^{*}}\frac{\bord f}{\bord n_{V_{\gamma^{*} } }    }(x)|\leq \sum_{x\in\gamma^{*}}|\frac{\bord f}{\bord n_{V_{\gamma^{*} } }   }(x)| .
\end{equation}

For every $x\in \gamma^{*}$ (viewed now as a vertex in $\delta
(V_{\gamma^{*}})$) we have that  $$\frac{\bord f}{\bord
n_{V_{\gamma^{*} }}  }(x)=(\vec{c}_{ \delta(V_{\gamma^*})  }(x), \overrightarrow{\frac{\bord f}{\bord n_{ V_{\gamma^{*} }     }}}(x)).$$ Hence by the
Cauchy-Schwartz inequality we have that 
\begin{equation}
\label{eq:CahchySchwartz}
| \frac{\bord f}{\bord n_{ V_{\gamma^{*} }     }}(x)
|=|(\vec{c}_{ \delta(V_{\gamma^*})  }(x), \overrightarrow{\frac{\bord f}{\bord n_{ V_{\gamma^{*} }     }}}(x))|\leq
\|\vec{c}_{ \delta(V_{\gamma^*})  }(x)\|\|\overrightarrow{\frac{\bord f}{\bord n_{ V_{\gamma^{*} }     }}}(x)\|.
\end{equation}

It is also clear that for every $x \in \gamma^{*}$ which is different from $x_0$ or $x_n$, we have that 

\begin{equation}
\label{eq:boundaryvsinterior}
 \|\vec{c}_{ \delta(V_{\gamma^{*}})  }(x)\|\leq \|\vec{c}_{F}(x)\|\ \  \mbox{\rm and }\ \
\|\overrightarrow{\frac{\bord f}{\bord n_{ V_{\gamma^{*} }     }}}(x)\|\leq \rho(x).
\end{equation}

If $x=x_0$ or $x=x_n$  we have that

\begin{equation}
\label{eq:moreboundaryvsinterior}
 \|\vec{c}_{ \delta(V_{\gamma^{*}})  }(x)\|\leq \|\vec{c}_{\delta(F)}(x)\|\ \  \mbox{\rm and }\ \
\|\overrightarrow{\frac{\bord f}{\bord n_{ V_{\gamma^{*} }     }}}(x)\|\leq \rho(x).
\end{equation}

Hence, by summing over all $x\in \gamma^{*}$, the definition of
$M$, Equations~(\ref{eq:properties1}), (\ref{eq:stokes}), (\ref{eq:CahchySchwartz}), (\ref{eq:boundaryvsinterior}) and 
(\ref{eq:moreboundaryvsinterior}) we have that

\begin{equation}
\label{eq:at last!}
 \frac{I(f)}{g M}\leq \sum_{x\in \gamma^{*}}\rho(x)={\mbox{\rm length}}_{\rho}(\gamma^*),
\end{equation} 
 which is first assertion of the
theorem.

\medskip
Let $\gamma$ be a shortest thick path connecting $P_2$ and $P_4$. 
It is clear that we may assume that $\gamma$ is simple. Let
 $x_0=\gamma \cap P_2 $ and let $x_n=\gamma\cap P_4$.
 
\myfigurepdf{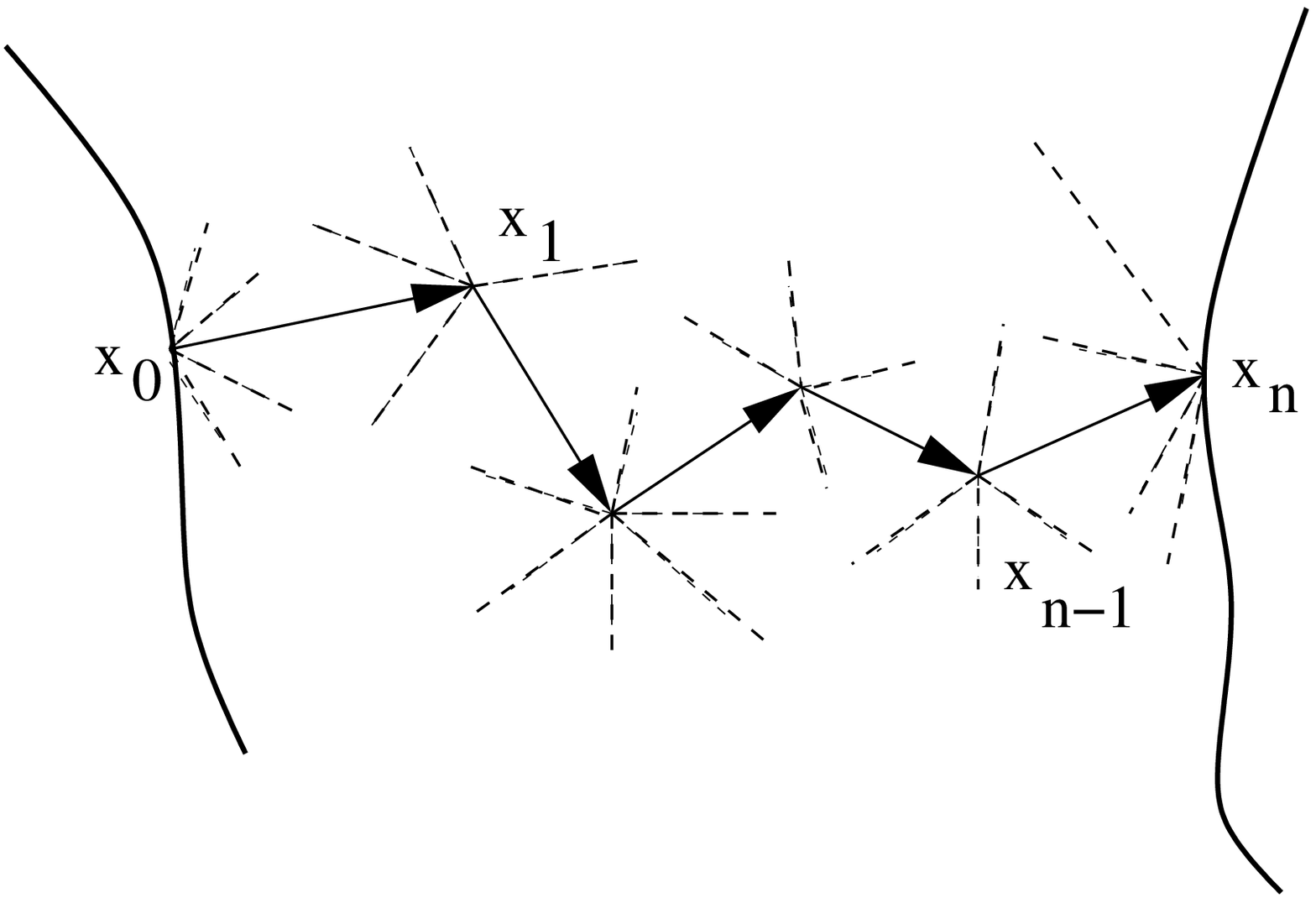}{6cm}{A shortest thick horizontal path between $P_2$ and $P_4$}

By integrating $\rho$ along $\gamma$ we have  that
$$\sum_{x\in\gamma}\rho(x)=\|\overrightarrow{\frac{\bord f}{\bord n_{F}}}(x_0)) \|+ \|\df(x_1)\|+\|\df(x_2)\|+\ldots+
\|\df(x_{n-1}\|+\| \overrightarrow{\frac{\bord f}{\bord n_F}}(x_n)) \|.$$

For each $i=1,\ldots, n-1$ we have that 
\begin{eqnarray}
\|\df(x_i)\| 
\geq 
\sqrt{c(x_i,x_{i+1})}|f(x_{i+1})-f(x_i) | \nonumber .
\end{eqnarray} 
It is easy to see that
 $ \|\overrightarrow{\frac{\bord f}{\bord n_{F}}}(x_0)) \| \geq \sqrt{c(x_0,x_1)}\,|f(x_1)-f(x_0)|$
and that $\|\overrightarrow{\frac{\bord f}{\bord n_{F}}}(x_n)) \|  \geq  \sqrt{c(x_{n-1},x_n)}\,|f(x_n)-f(x_{n-1})|$.  We now sum over all $x_i$, use
the definition of $m$ and the triangle inequality to obtain that

\begin{eqnarray}
\sum_{x\in\gamma}\rho(x) &\geq &\nonumber
\sum_{i=0}^{n-1}\sqrt{c(x_i,x_{i+1})}|f(x_{i+1})-f(x_i) |\\
&\geq & m |f(x_1)-f(x_0)+f(x_2)-f(x_1) +\ldots+
f(x_n)-f(x_{n-1})| \\ \nonumber& \geq &m g.
\end{eqnarray}

Assertion\ (2)\ of the theorem now follows.
\eop{\ref{le:length estimate}}

\rem It is easy to check that Assertion\ (2)\ of the theorem will hold for a larger class of horizontal paths. 

\medskip
We now provide an upper bound for the product of the lengths of any shortest paths in the network.

\blemm
\label{le:upper bound}
Let $\T$ be a triangulation of a topological quadrilateral. Let  $\bar{\Gamma}( F)$ be the associated network.
Let $f$ be the Dirichlet-Neumann boundary value function with some constant $g$.  Let $\rho$ be the gradient metric induced by $f$.  Then for any  $\rho$ shortest curves $\alpha,\beta$ in $\bar{\Gamma}(F)$ we have that 
 $$ {\mbox{\rm length}}_{\rho}(\alpha)    {\mbox{\rm length}}_{\rho}(\beta)  \leq l(|V|)  I(f) ,$$ where $l(|V|)$ is some constant which depends on $|V|$.
 \elemm
 \dem Let $\alpha=(x_0,x_1,\ldots,x_n)$ be a a $\rho$ shortest curve in $\bar{\Gamma}( F)$ connecting the vertex $x_0$ to the vertex $x_n$. Then
${\mbox{\rm length}}_{\rho}(\alpha)= \sum_{x\in\alpha}\rho(x)$. By the definition of $\rho$ (Definition~\ref{de:gradient metric})  we  have for all $x\in \bar F$ that

$$\rho(x) =\left(  \sum_{y\in \bar F} c(x,y)(f(x)-f(y))^2 \right)^{1/2} \leq 
\left( \sum_{(x,y)\in \bar E} c(x,y)(f(x)-f(y))^2 \right)^{1/2}= \sqrt{I(f)}.$$

It follows from Chapter 31 in \cite{Sed} (with only minor changes needed in our setting) that $n=O((|E|+|V|)\log|V|)$. Since $\bar{\Gamma}( F)$  is planar we also have that $|E|=O(|V|)$. Hence we have that $n=O(|V|\log|V|)$. Therefore it follows that $${\mbox{\rm length}}_{\rho}(\alpha)= O(|V|\log|V|)
\sqrt{I(f)}.$$ The assertion of the lemma follows easily.\hfill\eop{\ref{le:upper bound}}

\bcoro
\label{co:product of lengths}
Under the assumptions of the Theorem~\ref{le:length estimate} and Lemma~\ref{le:upper bound}  we have that 
$$l(|V|)  I(f)\geq {\mbox{\rm length}}_{\rho}(\gamma^*){\mbox{\rm length}}_{\rho}(\gamma)\geq \frac{m}{ M}I(f).$$
\ecoro
\rem
In the case that $c(x,y)\equiv 1$, it is easy to see that $\frac{m}{M}=\frac{1}{\sqrt{k}}$.

\section{An example}
\label{se:example}
\myfigurepdft{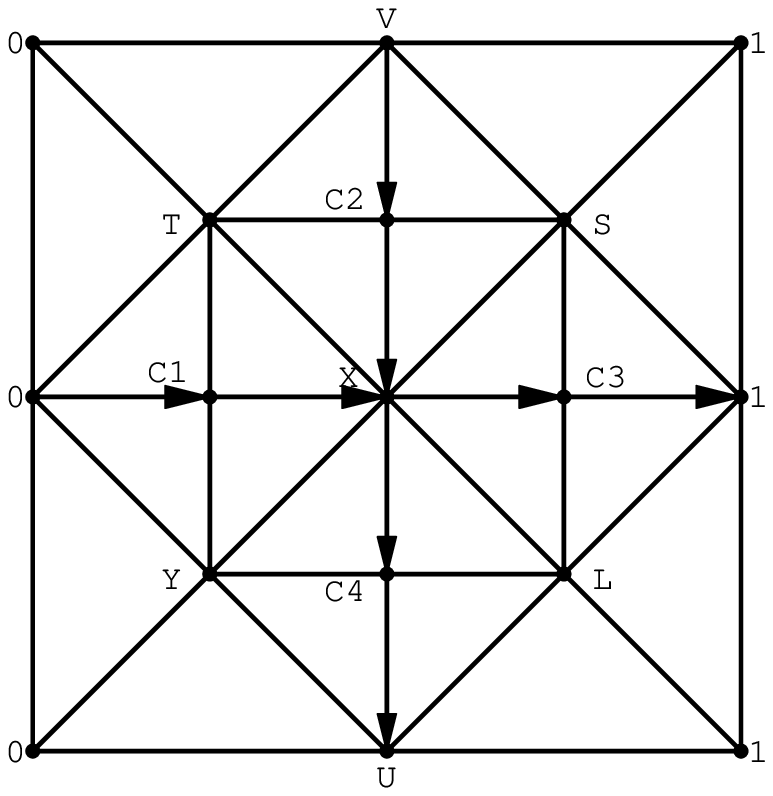}{6.5cm}{}
With the triangulation above, let us solve the Dirichlet-Neumann boundary value  (Definition~\ref{de:boundary function}) with $c(x,y)\equiv 1$ and $g=1$.
By abuse of notation let us use the same letter to indicate both the vertex name and the value of the solution at this vertex. Simple calculations (performed with Mathematica) shows the following.
\begin{enumerate} 
\item The solution is $(X,V,S,T,Y,L,U,C1,C2,C3,C4)=(\frac{1}{2},\frac{1}{2}, \frac{31}{44}, \frac{13}{44}, \frac{13}{44}, \frac{31}{44},\frac{1}{2}, \frac{3}{11},\frac{1}{2},\frac{8}{11},\frac{1}{2})$.
\item $[0,C1,X,C3,1]$ is the shortest thick horizontal geodesic and its length, denoted by $l_{h}$, with respect to the gradient metric is approximately $2.23111 $,
\item $[V,C2,X,C4,U]$ is the shortest thick vertical geodesic and its length, denoted by $l_v$, with respect to the gradient metric is approximately  $1.67733 $,
\item the energy, $I(f)$ of the solution equals $16/11\sim1.45455$,
\item $k=8$, and
\item  $\frac{ l_{h}  l_{v} }{I(f)}-\frac{1}{\sqrt{8}} \sim 2.21929 >0$.

\end{enumerate}

We conclude that the lower bound provided by the remark following Corollary~\ref{co:product of lengths} is not sharp.

\medskip
\begin{flushleft}
{\small 
University of Georgia\\
Dept. of Mathematics\\
Athens, GA 30602\\
saarh@math.uga.edu
}

\end{flushleft}

\end{document}